\documentclass{amsproc}

\copyrightinfo{2001}{American Mathematical Society}

\newtheorem{theorem}{Theorem}

\theoremstyle{definition}

\theoremstyle{remark}

\newcommand{\vv}{{\bf v}}

\newcommand{\what}{\widehat}

\newcommand{\R}{{\bf R}}
\newcommand{\Z}{{\bf Z}}
\newcommand{\Q}{{\bf Q}}

\newcommand{\q}{q}

\newcommand{\br}[1]{\| #1 \|}

\begin{document}

\title{Random walks with badly approximable numbers}

\author{Doug Hensley}
\address{Department of Mathematics, Texas A\&M University,
 College Station, TX  77843}
\email{Doug.Hensley@math.tamu.edu}

\author{Francis Edward Su}
\address{Department of Mathematics, Harvey Mudd College,
  Claremont, CA  91711}
\email{su@math.hmc.edu}
\thanks{Francis Su is grateful to the School of Operations Research at 
Cornell University for their hospitality during a sabbatical in which
this was completed.}


\begin{abstract}
Using the discrepancy metric, we analyze the rate of convergence
of a random walk on the circle generated by $d$ rotations, and
establish sharp rates that show that badly approximable $d$-tuples
in $\R^d$ give rise to walks with the fastest convergence.

\end{abstract}


\maketitle

\section{Introduction}

Fix $\bar\alpha=(\alpha_1, \alpha_2, ..., \alpha_d)$,
a $d$-tuple of real numbers, not all rational.
Consider a random walk on the circle
$S^1$ which proceeds as follows: at each step one of the $\alpha_i$
is chosen (with probability $1/d$) and the walk moves forwards or
backwards (with probability $1/2$)
on the circle by an angle $2\pi
\alpha_i$.  The distribution of this walk
converges to the uniform (Haar) measure on $S^1$; in this paper,
we investigate how quickly this convergence can take place.

To quantify this, identify $S^1$ with the unit circle in $\R^2$, and
let $Q$ be the
probability measure supported on $\{e^{\pm 2\pi i\alpha_j} \}$, $j=1,...,d$,
with equal weights $1/2d$.
If at time $k=0$ the random walk starts at $1$, then the convolution power
$Q^{*k}$ represents the probability distribution of the random
walk at time $k$.
In this framework, we study how quickly $Q^{*k}$ approaches
$U$, the uniform (Haar) measure on the unit circle.
Define
the {\em discrepancy} $D(P)$ of a probability measure on $S^1$ by
\begin{equation}
\label{discrepdefn}
D(P) = \sup_I |P(I) - U(I)|
\end{equation}
where $I$ represents any connected interval on $S^1$.  The
discrepancy $D(P)$ measures how uniform $P$ is, and metrizes 
weak-* convergence to the uniform distribution on $S^1$.
The usual definition of discrepancy found in the study of uniform
distribution of sequences mod 1 (see \cite{kuipers-nied}) is a
special case of our definition, by letting $P$ be an {\em equally
weighted} measure on a sequence.

We seek bounds on how quickly $D(Q^{*k})$ diminishes as a
function of the step size $k$ and the numbers $\alpha_i$.  Clearly the
convergence will not occur if all $\alpha_i$ are rational.

The case $d=1$ reduces to a random walk generated by one
irrational rotation; this problem was posed by Diaconis in
\cite[p. 34]{diaconis}. Su \cite{su} obtained discrepancy bounds
for this walk, showing that it converges most quickly when
$\alpha$ is a quadratic irrational. In that case, $$
\frac{C_1}{\sqrt{k}} \leq D(Q^{*k}) \leq \frac{C_2}{\sqrt{k}} $$
for constants $C_1, C_2$ which can be determined given $\alpha$.
We investigate whether the inclusion of additional generators can
speed up the convergence of this walk. In particular we show that
if $\bar\alpha$ is a {\em badly approximable} vector in $\R^d$,
then the random walk generated by the $\bar\alpha$ converges as
fast as possible for a set of $d$ generators. In Theorem
\ref{thm1} we obtain matching upper and lower bounds:
$$
\frac{C_1}{k^{d/2}} \leq D(Q^{*k}) \leq \frac{C_2}{k^{d/2}}.
$$
We also show how the constants depend on $\bar\alpha$.

This result is reminiscent of results in the theory of uniform
distribution of sequences mod 1; however, there is a notable lack
of any log terms in the upper and lower bounds.  To compare for
$d=1$, it is known (e.g., see \cite{kuipers-nied}) that for a
fixed irrational $\alpha \in \R$ the discrepancy of the sequence
$\{\alpha,2\alpha,3\alpha,\dots,k\alpha\}$ in $\R \mod 1$ falls as
$\log k/k$ up to constant factors.  For our corresponding random
walk, it is not a surprise that the exponent gets halved, but it
is quite surprising that the log term disappears.

While the total variation metric is more frequently used to study
random walk convergence (e.g., see \cite{diaconis}), it is not
appropriate for this walk, because $Q^{*k}$ is a finitely
supported measure and its total variation distance from $U$
remains at 1 for all $k$.  Moreover, we favor the use of
discrepancy over other common metrics on probabilities (e.g.,
Wasserstein, Prokhorov) because bounding techniques for such
metrics are not well developed, and the use of Fourier bounds for
the discrepancy metric admits the possibility of sharp analysis.
An analogous discrepancy metric may be used to study random walks on other
spaces, e.g., \cite{su-thesis, su-leveque, su-drunkard}.
For a survey of bounds relating various metrics, see
\cite{gibbs-su}.

We also remark that while the literature contains many results on
rates of convergence for random walks on finite groups (e.g., see
\cite{diaconis}), and a few results for walks on infinite compact
groups (\cite{porod},\cite{rosenthal}), very little has been done
for {\em discrete} walks on infinite compact groups, mainly due to
lack of bounding techniques for appropriate metrics.    Our
analysis of the $d$-generator random walk on the circle reveals
that bounds for the discrepancy metric are refined enough to yield
sharp rates of convergence; this and \cite{su} are the only sharp
results we are aware of in this direction.

\section{Random Walk Bounds}
Throughout this paper, let $\br{{\bf x}}$ denote the Euclidean
distance of ${\bf x} \in \R^d$ to the nearest integer lattice
point (the dimension inherent in the notation $\br{\cdot}$ is to
be understood by context). Given an arbitrary $d$-tuple
$\bar\alpha=(\alpha_1, \alpha_2, ..., \alpha_d)$, the Dirichlet
Approximation Theorem \cite[p.68]{drmotatichy} implies that there
exists a constant $B=B(\bar\alpha)$ such that for any $\q$, there
exists a positive integer $N \leq \q^d$ such that
\begin{equation}
\label{dirichletbound} \br{N \bar \alpha} < \frac{B}{\q} \leq
\frac{B}{N^{1/d}}.
\end{equation}
In fact for any $d$-tuple $\bar\alpha$, the above bound holds for
$B=\sqrt{d}$.

The $d$-tuple $\bar\alpha$ is said to be {\em badly approximable}
if there exists a constant $\beta=\beta(\bar\alpha)>0$ such that
for all positive integers $N$,
\begin{equation}
\label{thebound}
\br{N \bar\alpha}
\geq \frac{\beta}{N^{1/d}}.
\end{equation}

We refer to $B$ and $\beta$ as {\em approximation constants} for
$\bar\alpha$.  The approximation constant $\beta$ is only defined
for badly approximable $\bar\alpha$.

Our main result is the following:

\bigskip
\begin{theorem}
\label{thm1} Suppose $\bar\alpha=(\alpha_1, \alpha_2, ...,
\alpha_d)$ is badly approximable.  Let $Q$ denote the generating
measure of the random walk on $S^1$ generated by the $d$
generators $\alpha_i$. Then the discrepancy of the $k$-th step
probability distribution of the walk satisfies
\begin{equation}
\label{thmbound}
\frac{C_1}{k^{d/2}} \ \leq D(Q^{*k}) \leq \ \frac{C_2}{k^{d/2}}
\end{equation}
with constants that depend on $d$ and $\bar\alpha$:
\begin{eqnarray*}
C_1 &=& 0.0947 \  (\sqrt{d}/5B)^d \\
C_2 &=& 19.857 \ (d\sqrt{d}/\beta)^d
\end{eqnarray*}
where $B,\beta$ are approximation constants for $\bar\alpha$.
The lower bound holds for walks generated by an arbitrary $d$-tuple.
\end{theorem}

\bigskip
Therefore, among all $d$-tuples, walks generated by badly
approximable $d$-tuples converge the fastest.  Note that the
constants $C_1, C_2$ depend on the approximation constants of
$\bar\alpha$.  However since one may choose $B = \sqrt{d}$, it
follows that $C_1 \geq 0.0947 (1/5)^d$.

\begin{proof}
We establish the lower bound using an inequality of
Su \cite{su-leveque} for the discrepancy of an arbitrary
probability measure $P$ on $S^1$:
\begin{equation}
\label{su-lower-bound}
 D(P) \ \geq \ \left( \frac{2 }{\pi^2 }
\sum_{m=1}^{\infty} \frac{|\what{P} (m)|^2}{ m^2 } \right)^{1/2} .
\end{equation}
Here $\what{P}(m)$ denotes the $m$-th Fourier coefficient of $P$
(viewed on $\R$ mod $\Z$). Since every term is positive, we can
use a dominant term in the sum as a lower bound.

The Fourier coefficients for the $d$-generator walk are
\begin{equation}
\label{f-coeff}
 \what{Q}(m) = \sum_{l=1}^d \frac{1}{2d} (e^{2\pi
im\alpha_l} + e^{-2\pi im\alpha_l})
 = \frac{1}{d} \sum_{l=1}^d
\cos(2\pi m \alpha_l).
\end{equation}
Since
$\cos(2\pi x)=\cos(2\pi\br{x})\geq 1-2\pi^2 \br{x}^2$, for any $m$
we have
\begin{equation}
\label{circ-a}
\what{Q}(m)
\geq 1 - \frac{2\pi^2}{d} \br{m \bar \alpha}^2
\end{equation}
The inequality $(1-x)^k \geq 1-kx$ holds for $k\geq 1$ and $x\leq
1$. Then $$ |\what{Q}(m)|^k \geq  1 - \frac{2\pi^2 k}{d} \br{m
\bar \alpha}^2 $$ as long as $2\pi^2 k \br{m \bar \alpha}^2 /d <
1$. We ensure this by setting $Z_1<1$ (to be specified shortly)
and let $\q= (2 \pi^2 B^2 k / Z_1 d)^{1/2}$.  Then
(\ref{dirichletbound}) implies there exists an integer $m \leq
\q^d$ such that $\br{m \bar \alpha} < B/\q$, which yields $2\pi^2
k \br{m \bar \alpha}^2 /d < Z_1 < 1$, as desired.  (Note that $\q$
was chosen to ensure $|\what{Q}(m)|^k \geq 1-Z_1$.)

For this $m$ we use just the $m$-th term in
(\ref{su-lower-bound}), and since $m<\q^d$ and $\what{Q^{*k}}(m) =
\what{Q}^k(m)$, we obtain
\begin{equation}
\label{circ-b}
D(Q^{*k})
\geq \frac{\sqrt{2} |\what{Q}(m)|^k}{\pi m}
\geq C_1 \ k^{-d/2}
\end{equation}
with
$$C_1 = \frac{\sqrt{2}(1-Z_1)}{\pi} \left(\frac{Z_1
d}{2\pi^2 B^2}\right)^{d/2}. $$
Choosing $Z_1=2\pi^2/25 < 1$, we
recover the lower bound (\ref{thmbound}) of the theorem.

The upper bound is trickier to compute.
We start with the Erd\"os-Tur\'an inequality \cite{nied-phil}:
given a probability measure $P$ on $S^1$, for any integer $M$,
\begin{equation}
\label{erdos-turan}
D(P)  \leq  \frac{4}{M+1} + \frac{4}{\pi}
        \sum_{m=1}^{M} \frac{ | \what{P}(m) |}{m}
\end{equation}
where
$\what{P}$ represents the Fourier transform of $P$.
Note that one may choose $M$ in the Erd\"os-Tur\'an inequality so as
to optimize the bound obtained.

Since $| \cos(2\pi x)| \leq 1 - 4 \br{2 x}^2 $ for all $x \in \R$, it
follows from (\ref{f-coeff}) that
\begin{eqnarray*}
| \what{Q}(m) |
&\leq & \frac{1}{d} \sum_{l=1}^d |\cos(2\pi m \alpha_l)|\\
&\leq & 1 - \frac{4}{d} \sum_{l=1}^d \br{ 2 m \alpha_l }^2 \\
&\leq & \exp(- \frac{4}{d} \br{2 m \bar \alpha }^2 ).
\end{eqnarray*}

In light of the Erd\"os-Tur\'an inequality and
$\what{Q^{*k}}(m) = \what{Q}^k(m)$,
we need to estimate a sum of the form
$$
\sum_{m=1}^{M} \frac{ | \what{Q}^k(m) |}{m}
\leq \sum_{m=1}^{M}
\frac{1}{m} \exp(-\frac{4k}{d} \br{2 m \bar\alpha}^2) =: S.
$$

Since $M$ may be chosen freely, choose an integer $M$ such that
\begin{equation}
\label{M-defn}
M \leq \frac{1}{2} ( \beta^2 k / d^3 )^{d/2} < M+1.
\end{equation}
Recall that $\beta$ is an approximability constant for $\bar\alpha$
from (\ref{thebound}) and $k$ is the number of
steps in the walk.
The reason for this choice of $M$ will be evident later.
Choose an integer $J$ such that
\begin{equation}
\label{J-defn}
2^{J-1} \leq M \leq 2^J -1.
\end{equation}

The sum in $S$ may be grouped into $J$ cohorts of integers $m \in
[2^{j-1}, 2^j-1]$ for $j=1,...,J$. Within each such cohort, the
use of inequality (\ref{thebound}) yields $\br{2 m \bar\alpha }
\geq \beta /(2m)^{1/d} \geq \beta / 2^{(j+1)/d}$. This says that
the points of the sequence $\{ 2m\bar\alpha\}$ (mod 1) in the unit
cube in $\R^d$ are bounded away from the corners of the unit cube.
In fact, they are bounded away from each other, since if $m_1, m_2
\in [2^{j-1}, 2^j-1]$, then
\begin{equation}
\label{points-far}
\br{2(m_1-m_2)\bar\alpha} \geq \beta / 2^{(j+1)/d}
\end{equation}
as well.
Therefore if $s=\beta/\sqrt{d}\,2^{(j+1)/d}$, then
any box of side length $s$
can contain at most one of the multiples from the sequence
$\{2 m \bar\alpha \}$ (mod 1), $m \in [2^{j-1}, 2^j-1]$, with no such multiples
in boxes containing the origin.

So divide up the unit $d$-cube into disjoint boxes of side
length $s$ and sides parallel to the axes.
In the worst case, all $M$ points are distributed in the
boxes in the $2^d$ corners
of the unit cube nearest the origin.  The nearest points in such boxes
(under the $L_1$ metric) are at integral multiples of $s$ and
extend out in layers to distance at most
$M^{1/d}$.  A crude upper bound for the number of boxes whose
nearest point is at $L_1$-distance $ns$ from the origin is $(n+1)^{d-1}$,
and in the Euclidean metric the nearest point in such boxes is at
least distance $ns/\sqrt{d}$ from the origin.
Hence we can bound $S$ by grouping first by cohorts, then by corners
and layers:
\begin{eqnarray*}
S &\leq&
\sum_{j=1}^{J}
\sum_{m=2^{j-1}}^{2^j - 1}
   \frac{1}{m} \exp(-\frac{4k}{d} \br{2 m \bar\alpha}^2) \\
&\leq&
\sum_{j=1}^{J}  {2^d}
\sum_{n=1}^{M^{1/d}}
        \frac{(n+1)^{d-1}}{2^{j-1}}
        \exp(- \frac{4k}{d} \frac{\beta^2 n^2}{d \cdot 2^{2(j+1)/d}})\\
&\leq&
\sum_{j=1}^{J}  \frac{2^d}{2^{j-1}}
\sum_{n=1}^{\infty}
        (n+1)^{d-1}
        \exp(- 4 n^2 d 2^{2(J-j-1)/d}).
\end{eqnarray*}
The second inequality used the bound
$\br{2 m \bar\alpha} \geq ns/\sqrt{d}$ and the definition of $s$.
The third inequality follows by noting
$k \geq d^3 2^{2J/d} / \beta^2$
from the definitions of $J$ and $M$ in (\ref{M-defn}) and (\ref{J-defn}).

Using $j \leq J$, the log derivative of the
expression of the innermost sum with respect to
$n$ can be bounded:
$$\frac{d-1}{n+1} - 8 n d 2^{2(J-j-1)/d}
\leq \frac{d-1}{n+1} - 8 n d 2^{-2/d}
\leq \frac{d-1}{2} -  d
\leq -1
$$
for all choices of $n \geq 1$ and $d \geq 1$.
Hence the expression in the inner sum
decreases geometrically
(by at least ratio $e^{-1}$)
as $j$ increases, and so the inner sum is
bounded by the first term (at $n=1$)
times the constant $Z_2=1/(1-e^{-1}) \approx 1.5820$.
Thus
$$ S \leq \sum_{j=1}^{J}  Z_2 \,
        2^{2d-j}
        \exp(- 4 d 2^{2(J-j-1)/d} ).
$$
This sum may be bounded by noting that the largest term
occurs when $j=J$.
For $j\leq J$, the log derivative of the terms
with respect to $j$ is
$\ln 2(-1 + 8 \cdot 2^{2(J-j-1)/d})\geq \ln 2$ for $d\geq 1$.  Thus
the sum decreases geometrically (with at least ratio $1/2$) as $j\leq J$
decreases, so the sum is therefore bounded by twice the
final term at $j=J$:
\begin{eqnarray}
\nonumber
S
&\, \leq \, &
2 Z_2 \, 2^{2d-J} \exp(- 4 d 2^{-2/d})\\
\label{maxbound}
& \leq &
\frac{ Z_2 2^{2d+1} \exp(-4 d 2^{-2/d})}{M+1}
\leq \frac{4.6559}{M+1}
\end{eqnarray}
where the second inequality follows from (\ref{J-defn}), and the
final inequality from noting that the numerator is greatest for
$d=1$. Using the Erd\"os-Tur\'an inequality we obtain the upper
bound $$ D(Q^{*k}) \leq \frac{4}{M+1} + \frac{4}{\pi} S \leq
\frac{ 4 + (4/\pi) 4.6559}{M+1} \leq \frac{ 9.9281 }{M+1} $$ for
all $d\geq 1$. An application of (\ref{M-defn}) produces $$
D(Q^{*k}) \leq 19.857\,d^{3d/2} \beta^{-d} k^{-d/2} $$ which
establishes the constant $C_2$ in the statement of the theorem.
\end{proof}

We remark that this argument simplifies the proof given in
\cite{su} for the case $d=1$.  For specific $d$, the constants
$C_1, C_2$ can be improved by adjusting the derivations of the
constants $Z_1, Z_2$ and the bound for the last inequality in
(\ref{maxbound}).

\section{Choosing Generators}
Badly approximable $d$-tuples are plentiful. Cassels
\cite{cassels} shows that there uncountably many badly
approximable $d$-tuples in $\R^d$.  Moreover, results of
Schmidt \cite[pp. 53-59]{schmidt} imply that the Hausdorff
dimension of the set of badly approximable vectors is positive.
For some concrete examples, it can be shown (see, e.g., \cite[p.
68]{drmotatichy}) that if $1, \alpha_1, ..., \alpha_d$ are
linearly independent over $\Z$, and if the degree of the extension
$[\Q(\alpha_1,...,\alpha_d):\Q] = d+1$, then $\bar\alpha$ is badly
approximable.

We have shown that badly approximable $d$-tuples in $\R^d$ give
rise to random walks with the fastest convergence, established
sharp rates of convergence, and exhibited how the constants depend
on the approximation constants $\beta, B$ of the given $d$-tuple.
However, just among badly approximable $d$-tuples, can we say
which random walks converge ``the fastest''?  By Theorem
\ref{thm1}, this amounts to identifying $d$-tuples with the
largest possible approximation constant $\beta$.

For $d=1$, the number $\phi=\frac{\sqrt{5}-1}{2}$ satisfies
(\ref{thebound}) and is therefore badly approximable; in fact, it
is known to be the ``most'' badly approximable number in the sense
that its approximation constant $\beta$ is larger than for any
other number.

For $d=2$, we conjecture that the ``most'' badly approximable
$2$-tuple is the vector $\vv= (\gamma^{-2},\gamma^{-1})$, where
$\gamma\approx 1.3247$ is the unique real root of $x^3-x-1$. For
this $\vv$, we believe (based on heuristic arguments and numerical
evidence) that for any $\beta < 0.54850...$,
there is a sufficiently large $N$ such that $n>N$ implies that
\begin{equation}
\label{lowerbound}
\| n \vv \| > \beta\ n^{-1/2}.
\end{equation}
Furthermore, no other vector can have a much larger approximation
constant $\beta$, because the work of Davenport and Mahler
\cite{davenport-mahler} implies that if $\beta
>(2/\sqrt{23})^{1/2}\approx .64577...$, then there is no vector
$\vv$ in $\R^2$ for which (\ref{lowerbound}) holds for all $n$.

Thus the pair $(\gamma^{-2},\gamma^{-1})$ yields a random walk on
the circle with 2 generators whose convergence rate appears close
to fastest possible over all badly approximable pairs. Is it
fastest? For larger $d$ the question is also open.

\end{document}